\documentclass[journal]{IEEEtran}
\usepackage[utf8]{inputenc}

\title{Profitable Emissions-Reducing Energy Storage}
\author{Daniel~J.~Olsen,~\IEEEmembership{Member,~IEEE},
        Daniel~S.~Kirschen,~\IEEEmembership{Fellow,~IEEE}}

\usepackage{amsmath}
\usepackage{graphicx}
\graphicspath{{figures/}}
\usepackage{array}                      
\usepackage{cite}                       
\usepackage{eurosym}                    
\usepackage{algorithm}
\usepackage{algpseudocode}              
\usepackage{ifthen}
\usepackage{textcomp}                   
\usepackage[hyphens]{url}               

\usepackage{tikz}
\usepackage{textcomp}
\usepackage{hyperref}
\usepackage{lipsum}

\newcommand{\copyrighttext}{%
    \footnotesize \textcopyright 2019 IEEE. Personal use of this material is permitted. Permission from IEEE must be obtained for all other uses, in any current or future media, including reprinting/republishing this material for advertising or promotional purposes, creating new collective works, for resale or redistribution to servers or lists, or reuse of any copyrighted component of this work in other works.	DOI: \href{https://ieeexplore.ieee.org/document/8844848}{10.1109/TPWRS.2019.2942549}, IEEE Transactions on Power Systems.}

\newcommand{\copyrightnotice}{%
    \begin{tikzpicture}[remember picture,overlay]
    \node[anchor=south,yshift=8pt] at (current page.south) {\fbox{\parbox{\dimexpr\textwidth-\fboxsep-\fboxrule\relax}{\copyrighttext}}};
    \end{tikzpicture}%
}

\newenvironment{ldescription}[1]
  {\begin{list}{}%
   {\renewcommand\makelabel[1]{##1\hfill}%
   \settowidth\labelwidth{\makelabel{#1}}%
   \setlength\leftmargin{\labelwidth}
   \addtolength\leftmargin{\labelsep}}}
  {\end{list}}

\begin{document}

\maketitle
\copyrightnotice

\begin{abstract}
    While energy arbitrage from energy storage can lower power system operating costs, it can also increase greenhouse gas emissions. If power system operations are conducted with the constraint that energy storage operation must not increase emissions, how does this constraint affect energy storage investment decisions? Two bi-level energy storage investment problems are considered, representing `philanthropic' (profitability-constrained) and profit-maximizing storage investors (PhSI, PMSI). A MILP heuristic is developed to obtain good candidate solutions to these inherently MINLP bi-level problems.
    
    A case study is conducted on a 30\% renewable system, with sensitivity analyses on the price of storage and the price of carbon emissions. Regardless of the emissions-neutrality constraint, a PhSI installs significantly more energy storage than a PMSI, increasing system flexibility. The effect of the emissions-neutrality constraint in the absence of a carbon price is to reduce the quantity of storage purchased and reduce annual system emissions (\texttildelow 3\%), with only minor increases in overall cost (\texttildelow 0.1\%). In cases with a carbon price, storage does not tend to increase emissions and the emissions constraint does not tend to decrease storage investment. The emissions-neutrality constraint is seen to deliver similar emissions reductions even in a system with much higher renewable penetration (46\%).
\end{abstract}


\let\thefootnote\relax\footnote{D. J. Olsen and D. S. Kirschen are with the University of Washington Department of Electrical Engineering in Seattle, WA 98195 USA (e-mail: \{djolsen, kirschen\}@uw.edu).}

\section*{Nomenclature}

\subsection*{Sets and Indices}
\begin{ldescription}{$xxxxx$}
    \item [$A$] Set of representative days, indexed by $a$.
    \item [$B$] Set of transmission network buses, indexed by $b$.
    \item [$I$] Set of generating units, indexed by $i$.
    \item [$L$] Set of transmission lines, indexed by $l$.
    \item [$S$] Set of generator power output blocks, indexed by $s$.
    \item [$T$] Set of time intervals, indexed by $t$ or $\tau$.
\end{ldescription}

\subsection*{Parameters}
\begin{ldescription}{$xxxxx$}
    \item [$b_{i,s}$] Marginal cost of block $s$ of generator $i$ (\$/MWh).
    \item [$C_{i}^{\text{min}}$] Minimum cost of generator $i$ (\$/h).
    \item [$C_{i}^{\text{su}}$] Start-up cost of generator $i$ (\$).
    \item [$d_{b,t,a}$] Demand at bus $b$, time $t$, day $a$ (MW).
    \item [${E}^{\text{max}}$] Regulator's GHG emission target (tons).
    \item [$E_{i}^{\text{min}}$] Minimum GHG emissions of generator $i$ (tons/h).
    \item [$E_{i}^{\text{su}}$] Start-up GHG emissions of generator $i$ (tons).
    \item [$f_{l}^{\text{max}}$] Capacity of transmission line $l$ (MW).
    \item [$g_{i}^{\text{max}}$] Maximum power output of generator $i$ (MW).
    \item [$g_{i}^{\text{min}}$] Minimum power output of generator $i$ (MW).
    \item [$g_{i,s}^{\text{max}}$] Maximum power output of block $s$, generator $i$ (MW).
    \item [$g_{i}^{\text{down}}$] Minimum down-time of generator $i$ (h).
    \item [$g_{i}^{\text{up}}$] Minimum up-time of generator $i$ (h).
    \item [$h_{i,s}$] Marginal GHG emissions of block $s$, generator $i$ (tons/MWh).
    \item [$m_{l,b}^{\text{line}}$] Line connection map. ${m}_{lb}^{\text{line}} = 1$ if line $l$ starts at bus $b$, $= -1$ if line $l$ ends in bus $b$, $0$ otherwise.
    \item [$m_{i,b}^{\text{unit}}$] Unit map. ${m}_{i,b}^{\text{unit}} = 1$ if generator $i$ is located at bus $b$, $0$ otherwise.
    \item [$P^{\text{CO}_2}$] GHG emissions tax rate (\$/ton-$\text{CO}_{2}$e).
    \item [$P^{\text{load}}$] Load shed penalty (\$/MWh).
    \item [$P^{\text{ren}}$] Renewable generation shed penalty (\$/MWh).
    \item [$w_{b,t}$] Renewable generation at bus $b$, time $t$ (MW).
    \item [$x_{l}$] Reactance of line $l$ ($\Omega$).
    \item [$\pi_{a}$] Probability of day $a$.
\end{ldescription}

\subsection*{Primal Variables}
\begin{ldescription}{$xxxxx$}
    \item [$C^{\text{gen}}$] System operator's generation cost (\$).
    \item [$C^{\text{shed}}$] System operator's shed cost (\$).
    \item [$E_{a}$] GHG Emissions for day $a$ (tons).
    \item [$E^{\text{total}}$] Total GHG emissions (tons).
    \item [$f_{l,t,a}$] Power flow on line $l$, time $t$, day $a$ (MW).
    \item [$g_{i,t,a}$] Power output of generator $i$, time $t$, day $a$ (MW).
    \item [$g_{i,s,t,a}$] Power output of generator $i$, block $s$, time $t$, day $a$ (MW).
    \item [$s_{b,t,a}^{\text{load}}$] Load shed at bus $b$, time $t$, day $a$ (MWh).
    \item [$s_{b,t,a}^{\text{ren}}$] Renewable generation shed at bus $b$, time $t$, day $a$ (MWh).
    \item [$u_{i,t,a}$] Binary variable for the commitment status of generator $i$, time $t$, day $a$.
    \item [$v_{i,t,a}$] Binary variable for the start-up of generator $i$, time $t$, day $a$.
    \item [$z_{i,t,a}$] Binary variable for the shut-down of generator $i$, time $t$, day $a$.
    \item [$\theta_{b,t,a}$] Voltage phase angle of bus $b$, time $t$, day $a$ (rad).
\end{ldescription}

\subsection*{Dual Variables}
\begin{ldescription}{$xxxxx$}
    \item [$\underline{\delta}_{i,s,t}$] Dual variable for generator segment lower limit constraint.
    \item [$\overline{\delta}_{i,s,t}$] Dual variable for generator segment upper limit constraint.
    \item [$\underline{\phi}_{i,t}$] Dual variable for renewable shedding lower limit constraint.
    \item [$\overline{\phi}_{i,t}$] Dual variable for renewable shedding upper limit constraint.
    \item [$\kappa_{b,t}$] Dual variable for storage state-of-charge tracking constraint.
    \item [$\underline{\xi}_{b,t}$] Dual variable for storage state-of-charge lower limit constraint.
    \item [$\overline{\xi}_{b,t}$] Dual variable for storage state-of-charge upper limit constraint.
    \item [$\underline{\rho}_{b,t}^{\text{chg}}$] Dual variable for storage charging power lower limit constraint.
    \item [$\overline{\rho}_{b,t}^{\text{chg}}$] Dual variable for storage charging power upper limit constraint.
    \item [$\underline{\rho}_{b,t}^{\text{dis}}$] Dual variable for storage discharging power lower limit constraint.
    \item [$\overline{\rho}_{b,t}^{\text{dis}}$] Dual variable for storage discharging power upper limit constraint.
    \item [$\lambda_{b,t}$] Dual variable for power balance constraint.
    \item [$\underline{\gamma}_{l,t}$] Dual variable for line flow lower limit constraint.
    \item [$\overline{\gamma}_{l,t}$] Dual variable for line flow upper limit constraint.
    \item [$\beta_{l,t}$] Dual variable for DC power flow constraint.
    \item [$\alpha$] Dual variable for emissions constraint.
\end{ldescription}

\subsection*{Variable Sets}
\begin{ldescription}{$xxxxx$}
    \item [$\Omega^{\text{C}}$] Set of binary commitment variables.
    \item [$\Omega^{\text{D}}$] Set of dispatch variables.
    \item [$\Omega^{\lambda}$] Set of dual dispatch variables.
    \item [$\Omega^{\text{S}}$] Set of storage investment variables.
\end{ldescription}

\section{Introduction} \label{sec_intro}

Increasing the penetration of renewable energy is a popular solution to decarbonization of power systems; for example, many jurisdictions have renewable portfolio standards, and pricing of greenhouse gas (GHG) emissions incentivizes more installation of renewable generation. Renewable energy also has the benefits of low operating cost and decreased reliance on fuel imports which may come from politically volatile regions. However, the most abundant sources of cost-effective renewable energy--wind and solar photovoltaic--suffer from uncertainty and variability in their power production. Grid-scale energy storage is often seen as a promising solution for the intermittency of renewable resources, and therefore a valuable contribution towards broader decarbonization efforts (in combination with other approaches such as demand response).
Unfortunately, while storage \textit{can} be used to reduce the carbon intensity of power system operations, studies have shown that under current market structures and generation mixes the use of energy storage can increase overall GHG emissions. 

This effect is primarily caused by the economic incentives to charge storage at times when the marginal generator is coal-fired (cheap, but carbon-intensive) and discharge at times when the marginal generator is gas-fired (more expensive, but less carbon-intensive). Although there are also economic incentives to improve efficiency without fuel substitution, or to charge using renewable generation which would otherwise be spilled, the positive emissions impact of these incentives may be outweighed by the negative emissions impact of gas-to-coal fuel substitution. In the absence of a price on GHG emissions, estimates of the overall impact of energy storage participation include a range of 104-407 kg $\text{CO}_2$ per MWh delivered for grid-scale storage \cite{hittinger2015} and 75-270 kg/MWh for behind-the-meter storage \cite{fisher2017}. How to correct for these increases in emissions is a topic of active consideration by regulators \cite{cpuc2019}.

This paper investigates the impact of enforcing an emissions-neutrality constraint (\textit{i.e.} system emission with storage participation shall not be greater than they would be without storage) on the operation of power systems, investment in grid-scale energy storage, and resulting power system emissions.

The modeling of energy storage operation and its impact on grid emissions has been studied using a wide variety of power systems models and modes of energy storage participation; consequently, findings on emissions impacts are varied as well.

Energy storage has been modeled as providing energy \cite{ummels2008, sioshansi2011, carson2013, babacan2018, zheng2018, hittinger2015, peterson2011, hoehne2016, arciniegas2018}, reserves \cite{noori2016, ryan2018, lin2016}, or both energy and reserves \cite{swider2007, das2015, sioshansi2009, fisher2017, craig2018, tuohy2011}.
Storage providing energy has been shown to be capable of increasing emissions \cite{ummels2008, sioshansi2011, carson2013, babacan2018}, or can either increase or decrease emissions depending on charge-scheduling heuristics \cite{hoehne2016, peterson2011, zheng2018} and generation mixture \cite{hoehne2016, peterson2011, hittinger2015, arciniegas2018}. In \cite{sioshansi2011}, higher penetrations of renewables amplify this emission-increasing effect.
Storage providing reserves has been shown to be capable of increasing emissions \cite{ryan2018}, decreasing emissions \cite{noori2016}, or can either increase or decrease emissions depending on reserve quantity, storage quantity, and reserve scheduling rules \cite{lin2016}.
Storage providing both energy and reserves has been shown to be capable of increasing emissions \cite{tuohy2011, fisher2017}, decreasing emissions \cite{swider2007, sioshansi2009, das2015}, or can either increase or decrease emissions depending on changing generation mixes \cite{craig2018}. In \cite{tuohy2011}, emissions are increased even in the presence of a 30 \euro{}/ton price on GHG emissions, although increasing wind penetration mitigates this effect.

When modeling storage's participation in power systems, many studies have treated the system's marginal emissions rate as fixed, and determined exogenously \cite{hoehne2016, noori2016, carson2013, hittinger2015, fisher2017, arciniegas2018, babacan2018, zheng2018}. This may be a reasonable approximation for small-scale energy storage, but larger quantities of energy storage will have the ability to change the marginal unit(s) of generation. More sophisticated approaches use an economic dispatch model \cite{peterson2011, sioshansi2011, lin2016} or a unit commitment \cite{sioshansi2009, tuohy2011, das2015, craig2018, ryan2018, swider2007, ummels2008}. Inclusions of the transmission network as in \cite{das2015, lin2016, ryan2018, swider2007, ummels2008} is also significant, as locationality affects storage operations in a transmission-constrained system.

Overall, the literature has mixed results on the overall emissions impact of energy storage. Identified factors which may influence whether energy storage is beneficial or detrimental from an emissions perspective include: charging/discharging strategies \cite{hoehne2016, peterson2011, zheng2018}, generation mixture \cite{hoehne2016, peterson2011, hittinger2015, arciniegas2018, craig2018}, reserve scheduling rules, \cite{lin2016}, and the quantity of storage installed \cite{lin2016}. Beyond the relevant `system' parameters, the direction of the emissions impact may also be determined by the particular day(s) under consideration, since the shape of the net load profile will determine the economic opportunities available to energy storage operators.

Previous work has investigated the emissions impact of storage participating in a market environment, but hasn't looked at how a daily emissions-neutrality constraint can impact storage investment decisions and resulting emissions. Lin \textit{et al.} \cite{lin2016} introduces an emissions-neutrality constraint for storage providing reserves, but only considers single-period economic dispatch, and do not assess the impact on storage profitability. Several papers have investigated storage expansion planning for low-carbon emissions goals \cite{haller2012, fursch2013, mileva2016}, but they do not address the emissions impact of adding storage to current market environments or storage profitability. A comprehensive review of storage expansion planning for low-carbon power systems is given in \cite{haas2017}.

By contrast, this paper investigates the impact of an emissions-neutrality constraint (ENC) on the profitability of various quantities of storage, and therefore the optimal quantity of storage to invest in, and finally the resulting GHG emissions impact.

Although prices on carbon emissions are typically considered to be the most efficient mechanisms to reduce emission quantity, in many places they are politically infeasible, or are set at quantities which are considered too low to achieve substantive reductions in GHG emissions \cite{worldbank2019}. Despite their relative inefficiency, alternative policy measures to reduce GHG emissions are often implemented instead.

Profitability is assessed from two merchant storage perspectives: a profit-maximizing storage investor (PMSI) or a `philanthropic' investor (PhSI) who only requires a specified minimum return on investment to cover installation costs. An example of an entity which may want to participate as a philanthropic storage investor would be a government with an interest in reducing costs and/or emissions from the power system in its jurisdiction while still participating in a competitive energy market with existing generation companies. Alternatively, a philanthropic emissions-neutral storage investor may be a public-private partnership; a governmental organization may offer low-cost financing, tax incentives, or other cost-reducing measure to a private storage investor in exchange for a commitment to offer storage dispatch control to the system operator with an ENC. These merchant energy storage perspectives are compared with a vertically-integrated utility (VIU) perspective, concerned only with overall costs and not storage profitability.

Though investments in energy storage by a system operator or a PMSI have reciprocal impacts on profitability \cite{pandzic2018}, in this framework we look at just one or the other in order to more clearly see the impact of the ENC. Although storage may earn revenue by participating in reserves, including reserves in power systems modeling introduces sensitivity to required reserves quantity and scheduling rules \cite{lin2016}. Therefore, reserves are omitted in order to focus on how an ENC impacts a system where energy participation can increase emissions.

The specific contributions of this paper are:

\begin{itemize}
    \item The formulation of two bi-level models (`philanthropic' and `profit-maximizing') to optimize merchant storage investments in light of emissions-neutrality constraints in commitment and dispatch.
    \item Development of a heuristic to quickly obtain good feasible solutions to this inherently non-convex and computationally difficult problem.
    \item Analysis of solution quality of the `philanthropic' problem using a relaxation that allows the evaluation of solution quality with reduced computational burden.
    \item Demonstration of these methods on a detailed case study.
    \item Sensitivity analysis showing the effect of governmental incentives or taxes on the optimal quantity of energy storage and the resulting operational impacts.
\end{itemize}

\section{Power System Model} \label{sec_power_sys_model}

A transmission-constrained unit commitment formulation is used to model a system operator's choice of online generators, the dispatch quantities of generators and energy storage, and the resulting prices.

The model implicitly assumes an adequately competitive electricity market, in which a centralized operator clears bids from generators and storage operators that accurately represent their true capacities and costs. The results in an adequately competitive bilateral market should be broadly similar, with the structural difference that the emissions-neutrality constraint would be enforced on the storage operator's purchases and sales. This system operation formulation is embedded into a multi-level optimization model used by an energy storage investor to decide the quantities and locations to install storage.

\subsection{System Operational Constraints}

Operational constraints for the system operator's unit commitment problem are given in \eqref{binary_2}-\eqref{emissions_constraint}. The constraints in this section apply for each representative day $a \in A$, however the index $a$ is omitted for brevity. Dual variables for each dispatch constraint are given in parentheses to the right of the equation.

\begin{gather}
    v_{i,t} - z_{i,t} = u_{i,t} - u_{i,t-1} \hspace{2pt} ; \hspace{2pt} \forall i \in I, t \in T
    \label{binary_2} \displaybreak[0] \\
    \sum_{\tau=t-\text{g}_{i}^{\text{up}}+1}^{t} v_{i,\tau} \le u_{i,t} \hspace{2pt} ; \hspace{2pt} \forall t \in T, i \in I
    \label{updown_1} \displaybreak[0] \\[-2pt]
    \sum_{\tau=t-\text{g}_{i}^{\text{down}}+1}^{t} z_{i,\tau} \le 1 - u_{i,t} \hspace{2pt} ; \hspace{2pt} \forall t \in T, i \in I
    \label{updown_2} \displaybreak[0] \\
    0 \le g_{i,s,t} \le g_{i,s}^{\text{max}} u_{i,t}
    \hspace{8pt} \forall i \in I, s \in S, t \in T \hspace{8pt} \left ( \underline{\delta}, \overline{\delta} \right )
    \label{seg_limit} \displaybreak[0] \\
    0 \le s_{b,t}^{\text{ren}} \le w_{b,t}
    \hspace{8pt} \left ( \underline{\phi}, \overline{\phi} \right ) \label{spill_limits} \\
    Q_{b,t} = Q_{b,t-1} + \eta J_{b,t}^{\text{chg}} - \frac{1}{\eta} J_{b,t}^{\text{dis}} \hspace{8pt} \forall b \in B, t \in T \hspace{8pt} \left ( \kappa \right )
    \label{soc_tracking} \displaybreak[0] \\
    0 \le Q_{b,t} \le Q_{b}^{\text{max}}
    \hspace{8pt} \forall b \in B, t \in T \hspace{8pt} \left ( \underline{\xi}, \overline{\xi} \right )
    \label{soc_limits} \displaybreak[0] \\
    0 \le J_{b,t}^{\text{chg}} \le J_{b}^{\text{max}}
    \hspace{8pt} \forall b \in B, t \in T
    \hspace{8pt} \left ( \underline{\rho}^{\text{chg}}, \overline{\rho}^{\text{chg}} \right )
    \label{charge_limits} \displaybreak[0] \\
    0 \le J_{b,t}^{\text{dis}} \le J_{b}^{\text{max}}
    \hspace{8pt} \forall b \in B, t \in T
    \hspace{8pt} \left ( \underline{\rho}^{\text{dis}}, \overline{\rho}^{\text{dis}} \right )
    \label{discharge_limits} \displaybreak[0] \\
    \sum_{i \in I} m_{i,b}^{\text{unit}} g_{i,t} + J_{b,t}^{\text{dis}} - \sum_{l \in L} m_{l,b}^{\text{line}} f_{l,t} + w_{b,t} - s_{b,t}^{\text{ren}} = \hspace{20pt} \nonumber \displaybreak[0] \\
    \hspace{40pt} d_{b,t} + J_{b,t}^{\text{chg}} \hspace{2pt} ; \hspace{2pt} \forall b \in B, t \in T
    \hspace{8pt} \left ( \lambda \right )
    \label{power_balance} \displaybreak[0] \\
    f_{l,t} = y_{l} \sum_{b \in B} m_{l,b}^{\text{line}} \theta_{b,t}  \hspace{2pt} ; \hspace{2pt} \forall l \in L, t \in T
    \hspace{8pt} \left ( \beta \right )
    \label{flow_def} \displaybreak[0] \\
    - f_{l}^{\text{max}} \le f_{l,t} \le f_{l}^{\text{max}} \hspace{2pt} ; \hspace{2pt} \forall l \in L, t \in T
    \hspace{8pt} \left ( \underline{\gamma}, \overline{\gamma} \right )
    \label{line_flow_limits} \displaybreak[0] \\
    E^{\text{total}} \le \chi E^{\text{baseline}} \hspace{8pt} \left ( \alpha \right )
    \label{emissions_constraint} \\
    g_{i,t} := g_{i}^{\text{min}} u_{i,t} + \sum_{s \in S} g_{i,s,t} \hspace{2pt} ; \hspace{2pt} \forall i \in I, t \in T 
    \label{gen_sum} \nonumber \\
    E^{\text{total}} := \sum_{i \in I} \sum_{t \in T} 
    \Big( E_{i}^{\text{min}} u_{i,t} + {E}_{i}^{\text{su}} v_{i,t}
    + \sum_{s \in S} h_{i,s} g_{i,s,t} \Big)
    \label{total_emissions_def} \nonumber
\end{gather}

Eqs. \eqref{binary_2}-\eqref{updown_2} constrain the generator commitment variables. Eq. \eqref{binary_2} relates status, startup, and shutdown variables ($u_{i,t}$, $v_{i,t}$, $z_{i,t}$, respectively). Eq. \eqref{updown_1} ensures that minimum up-times are respected, while \eqref{updown_2} ensures that minimum down-times are respected.
Generator heat rate curves are represented by piecewise-linear segments. Eq. \eqref{seg_limit} constrains the power in each segment $g_{i,s,t}$ to its segment capacity limit $g_{i,s}^{\text{max}}$.
Renewable spillage $s_{b,t}^{\text{ren}}$ is constrained in \eqref{spill_limits} to be non-negative and not more than the maximum renewable generation available $w_{b,t}$.
Energy storage charging $J_{b,t}^{\text{chg}}$ and discharging $J_{b,t}^{\text{dis}}$ are constrained in \eqref{soc_tracking}-\eqref{discharge_limits}, where \eqref{soc_tracking} tracks the state-of-charge of energy storage $Q_{b,t}$, \eqref{soc_limits} constrains the state-of-charge based on energy capacity $Q_{b}^{\text{max}}$, and \eqref{charge_limits}-\eqref{discharge_limits} constrain charging and discharging based on power capacity $J_{b}^{\text{max}}$. Storage charge and discharge efficiency are given by $\eta$.
The transmission network is represented by the DC power flow model in \eqref{power_balance}-\eqref{line_flow_limits}. Network topology is defined by matrices $m_{i,b}^{\text{unit}}$ for generators and $m_{l,b}^{\text{line}}$ for transmission lines. Eq. \eqref{power_balance} ensures power balance at each node in the network, \eqref{flow_def} relates flows $f_{l,t}$ to bus angles $\theta_{b,t}$ via line admittances $y_{l}$, and \eqref{line_flow_limits} constrains each line's power flow to its maximum magnitude $f_{l}^{\text{max}}$.
Finally, \eqref{emissions_constraint} represents the emissions-neutrality constraint (ENC), ensuring that GHG emissions $E^{\text{total}}$ (based on generator emissions at minimum power $E_{i}^{\text{min}}$, startup emissions $E_{i}^{\text{su}}$, and marginal emissions rates $h_{i,s}$) are not increased, relative to the emissions of the baseline (no-storage) solution $E^{\text{baseline}}$. These emissions rates may already be measured by existing air-quality regulations, and therefore available to the system operator \cite{silverevans2012}. To solve with a non-binding ENC, the value of $\chi$ is set to a large constant.

\subsection{System Operator's Problem}

The goal of the system operator is to minimize the total cost of supplying the demand of the system $C^{\text{gen}}$ (based on generator cost at minimum power $C_{i}^{\text{min}}$, startup costs $C_{i}^{\text{su}}$, and marginal costs $b_{i,s}$), subject to the operational constraints; this problem is formalized in \eqref{sysop_obj}-\eqref{sysop_constraints}.

\begin{gather}
    \min_{\Omega^{\text{C}}, \Omega^{\text{D}}} C^{\text{gen}}
    \label{sysop_obj} \\
    C^{\text{gen}} := \sum_{t \in T} \sum_{i \in I}
    \Big( C_{i}^{\text{min}} u_{i,t} + {C}_{i}^{\text{su}} v_{i,t} + \sum_{s \in S} b_{i,s} g_{i,s,t} \Big)
    \nonumber \\
    \text{subject to:} \nonumber \\
    \text{Eqs. \eqref{binary_2}-\eqref{emissions_constraint}}
    \label{sysop_constraints}
\end{gather}

\section{Investment Models} \label{sec_investment models}

The question of how much storage is `optimal' to install depends on the perspective taken. The simplest case is for a vertically-integrated utility (VIU), whose only objective is minimizing overall social cost. This perspective is simplest as decisions on storage investment $\Omega^{\text{S}}$, generator commitment $\Omega^{\text{C}}$, and dispatch $\Omega^{\text{D}}$ are conducted simultaneously. This perspective is formalized in \eqref{viu_obj}-\eqref{viu_batt_cost}.

The parameter $\pi_a$ represents the frequency of each representative day, while storage investment costs are determined by the amoritized per-MWh storage energy cost $c^{\text{Q}}$ and per-MW storage power cost $c^{\text{J}}$ for a given storage quantity. Although real costs for energy storage projects are not generally linear, a sub-linear relationship would be inherently non-convex and therefore complicate the optimization process. A linear model is consistent with projects with fixed project costs plus linear costs for capacity, with solutions differing only for small quantities of storage where the capacity costs do not dominate.

\begin{gather}
    \min_{\Omega^{\text{C}}, \Omega^{\text{D}}, \Omega^{\text{S}}} C^{\text{batt}} + \sum_{a \in A} \pi_{a} C_{a}^{\text{gen}}
    \label{viu_obj} \\
    \text{subject to:} \nonumber \\
    \label{viu_cons} \text{Eqs. \eqref{binary_2}-\eqref{emissions_constraint}} \\
    C^{\text{batt}} := \sum_{b \in B} \left ( c^{\text{Q}} Q_b^{\text{max}} + c^{\text{J}} J_b^{\text{max}} \right )
    \label{viu_batt_cost}
\end{gather}

This is contrasted with a centralized power market structure, in which a storage investor earns revenue based on locational marginal prices (LMPs), obtained from the commitment and dispatch solution determined by the system operator. A storage owner's net profit is determined by LMPs $\lambda_{b,t,a}$, power dispatch, and storage investment costs $C^{\text{batt}}$, as in \eqref{profit_def}.

\begin{equation}
    Profit := \sum_{a \in A} \left [ \pi_a \sum_{b \in B} \sum_{t \in T} \lambda_{b,t,a} \left ( J_{b,t,a}^{\text{dis}} - J_{b,t,a}^{\text{chg}} \right ) \right ] - C^{\text{batt}}
    \label{profit_def}
\end{equation}

A purely self-interested storage investor would have the sole objective of choosing storage investments to maximize their net profit, while LMPs are determined by a lower-level system operator determining unit commitment and dispatch. This perspective is formalized in \eqref{strat_storage_maxprofit_obj}-\eqref{strat_storage_maxprofit_ll}, and is referred to as the profit-maximizing storage investor (PMSI). Although owners of large geographically-dispersed storage installations may coordinate bids to maximize their energy-market profit \cite{mohsenian-rad2016}, the PMSI formulation assumes that once the storage investments are made, the energy market is competitive enough to be modeled with a cost-minimization problem.

\begin{gather}
    \max_{\Omega^{\text{S}}} Profit
    \label{strat_storage_maxprofit_obj} \\
    \text{subject to:} \nonumber \\
    \boldsymbol{\lambda} \in \text{arg } \min_{\Omega^{\text{C}}, \Omega^{\text{D}}} \Big \{ C^{\text{gen}}
    ; \text{subject to: Eqs. \eqref{binary_2}-\eqref{emissions_constraint}} \Big \}
    \label{strat_storage_maxprofit_ll}
\end{gather}

\begin{figure}[h]
    \centering
    \includegraphics[width=\linewidth]{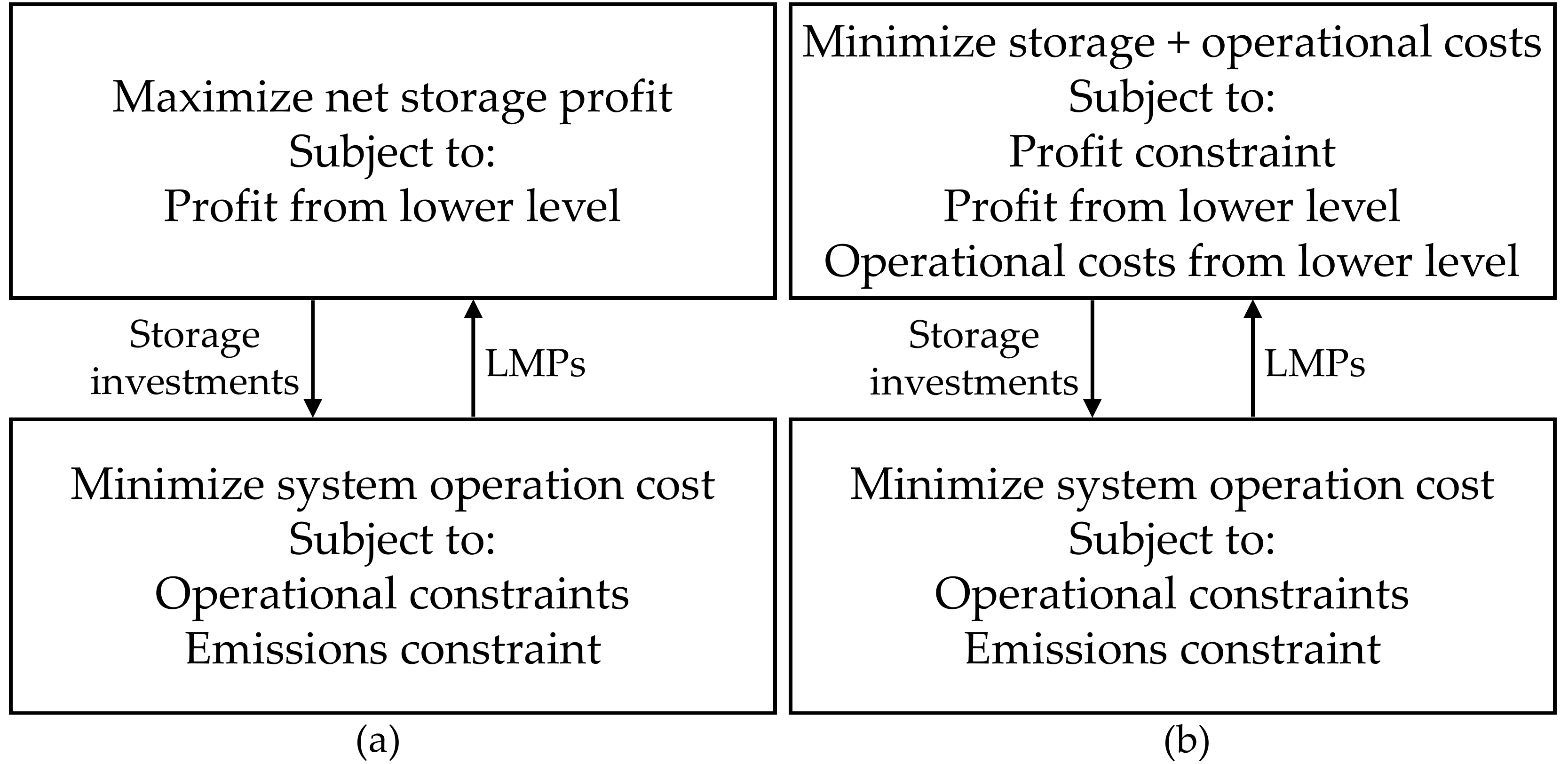}
    \caption{Bi-level formulations of merchant energy storage: (a) a profit-maximizing storage investor, (b) a `philanthropic' storage investor.}
    \label{strategic_storage}
\end{figure}

By contrast, there may be a storage investor who is less concerned with maximizing net profit than with lowering the overall social cost, subject to the constraint that enough energy market revenue is collected to cover the annualized investment cost \cite{dvorkin2017}. We refer to this perspective as a `philanthropic' storage investor (PhSI). Entities which may take this perspective include governmental or not-for-profit entities participating in a competitive energy market. This perspective is formalized in \eqref{strat_storage_phil_obj}-\eqref{strat_storage_phil_ll}. Both the PMSI and PhSI are bi-level optimization problems, and are illustrated in Fig. \ref{strategic_storage}.

\begin{gather}
    \min_{\Omega^{\text{S}}} C^{\text{batt}} + \sum_{a \in A} \pi_{a} C_{a}^{\text{gen}}
    \label{strat_storage_phil_obj} \\
    \text{subject to:} \nonumber \\
    Profit \ge 0
    \label{strat_storage_phil_profit_cons} \\
    C^{\text{gen}}, \boldsymbol{\lambda} \in \text{arg } \min_{\Omega^{\text{C}}, \Omega^{\text{D}}} \Big \{ C^{\text{gen}}
    ; \text{subject to: Eqs. \eqref{binary_2}-\eqref{emissions_constraint}} \Big \}
    \label{strat_storage_phil_ll}
\end{gather}

\section{Solution Methods} \label{sec_solution_method}

Solving the VIU perspective is a relatively straightforward MILP problem, since all decisions are made in a single level. However, the inclusion of energy market profit in the bi-level models complicates the solution process, since storage investment decisions influence energy market prices and energy market prices influence storage investment decisions. Bi-level optimization problems are inherently non-convex and NP-hard, even with a linear lower level \cite{dempe2015}; problems where the lower level is non-convex (\textit{e.g.} unit commitment) are even more challenging. For the purpose of this paper, we develop an iterative heuristic (Alg. \ref{alg_1}) in order to determine good candidate solutions to both bi-level problems in a reasonable amount of time. This algorithm is repeated for each combination of carbon price and effective storage price (described further in Section \ref{sec_case_study}). The quality of these solutions is evaluated in Section \ref{sec_solution_quality}.

All storage quantities $q$ between zero and the optimal VIU quantity $q^{\text{max}}$ are evaluated, and the resulting net storage profit and social cost for each quantity are evaluated once the unit commitment is solved. Finally, the profitable storage quantity with lowest social cost is chosen as the best candidate for the PhSI case, and the storage quantity with greatest net profit is chosen as the best candidate for the PMSI case.

\begin{algorithm}[H]
    \caption{Heuristic solution algorithm}
    \label{alg_1}
    \begin{algorithmic}[1]
        \Statex
        \State \Call{Define}{VIU model, Lower-level model}
        \State $q^{\text{max}} \gets \infty$ \Comment{$q$ is the system's total storage quanta}
        \State \Call{Solve}{VIU model}
        \State \Call{Record}{$q$, NetStorageProfit, SocialCost}
        \State $q^{\text{max}} \gets (q - 1)$
        \For{$i \gets 1 \text{\textbf{ to }} q^{\text{max}}$}
            \State \Call{Fix}{$q = i$}
            \State \Call{Solve}{Lower-level model}
            \State \Call{Record}{$q$, NetStorageProfit, SocialCost}
        \EndFor
        \State $q^{\text{philanthropic}} \gets $\Call{Select}{$q$}
        \Statex $\hspace{15pt} q \in \{ \text{argmin SocialCost}(q), \text{NetStorageProfit}(q) \ge 0 \}$
        \State $q^{\text{profit-max}} \gets $\Call{Select}{$q \text{ }|\text{ }q \in \text{argmax NetStorageProfit}(q)$}
    \end{algorithmic}
\end{algorithm}

\subsection{Assessing Solution Quality} \label{sec_solution_quality}

As the PhSI and PMSI problems are inherently non-linear and non-convex, evaluating the quality of candidate solutions (relative to a global lower bound) is not straightforward. However, there are two MILP relaxations of the PhSI problem, which can be used to find an upper bound for the optimality gap of candidate solutions. The first is the simply the VIU, since the profit constraint and the lower-level optimality constraint are both relaxed. The second relaxes the lower-level optimality constraint on the unit commitment variables to create a formulation which can be transformed into a profit-constrained single-level equivalent (PCSLE). The derivation of the PCSLE follows the approach of \cite{dvorkin2017}, and is summarized in Appendix A.

For the PMSI problem, obtaining a good estimate for the optimality gap is challenging. The same approach of relaxing the lower-level optimality constraint on the unit commitment variables would also yield a formulation which can be transformed into a single-level equivalent, but the relaxation would be significantly looser. This is because the upper-level objective is to maximize storage profit, so the relaxation will tend to choose uneconomic commitments which would never be chosen by a system operator, and therefore the upper bound on profit found by the relaxation will be very high. Instead, we evaluate local optimality by iteratively perturbing the heurtistic solution by one unit of storage at each bus, and ensuring that none of these alternative storage solutions produces a higher net profit.

\subsection{Limitations of the PCSLE} \label{sol_method_milp_discussion}

Although the MILP approximation of the original bi-level problem provides a relaxation which can assess the quality of the heuristic-found solutions, this formulation is not without its drawbacks: first and foremost is computational tractability. As mentioned in Section \ref{sec_solution_method}, single-level equivalents of bi-level problems are inherently NP-hard. In practice, solving problems with big-M approximations can be challenging, as the values of $M$ must be large enough to capture the full range of the continuous variables, but values which are too large cause difficulties for MILP solvers.

Second, as the MILP formulation is a relaxation, it may include integral solutions which are cheaper than the heuristic-found solution, but are not valid solutions to the original bi-level problem. A storage quantity which is not found to be profitable in the original bi-level problem may be profitable in this single-level relaxation, since unit commitment variables are no longer constrained to be part of the minimum-cost solution of the original lower-level problem. In essence, the single-level equivalent assumes that all decisions (storage quantity, unit commitments) are made simultaneously, and so a storage quantity which would lower total costs but otherwise be unprofitable can be made profitable by selection of a sub-optimal unit commitment solution. Therefore, the true optimum of the original bi-level problem may be greater than than the best lower-bound found by the single-level equivalent, and therefore the optimality gap for a given heuristic-found solution can only overestimate the true optimality gap.

\section{Case Study} \label{sec_case_study}

\begin{figure*}[h]
    \centering
    \includegraphics[width=\linewidth]{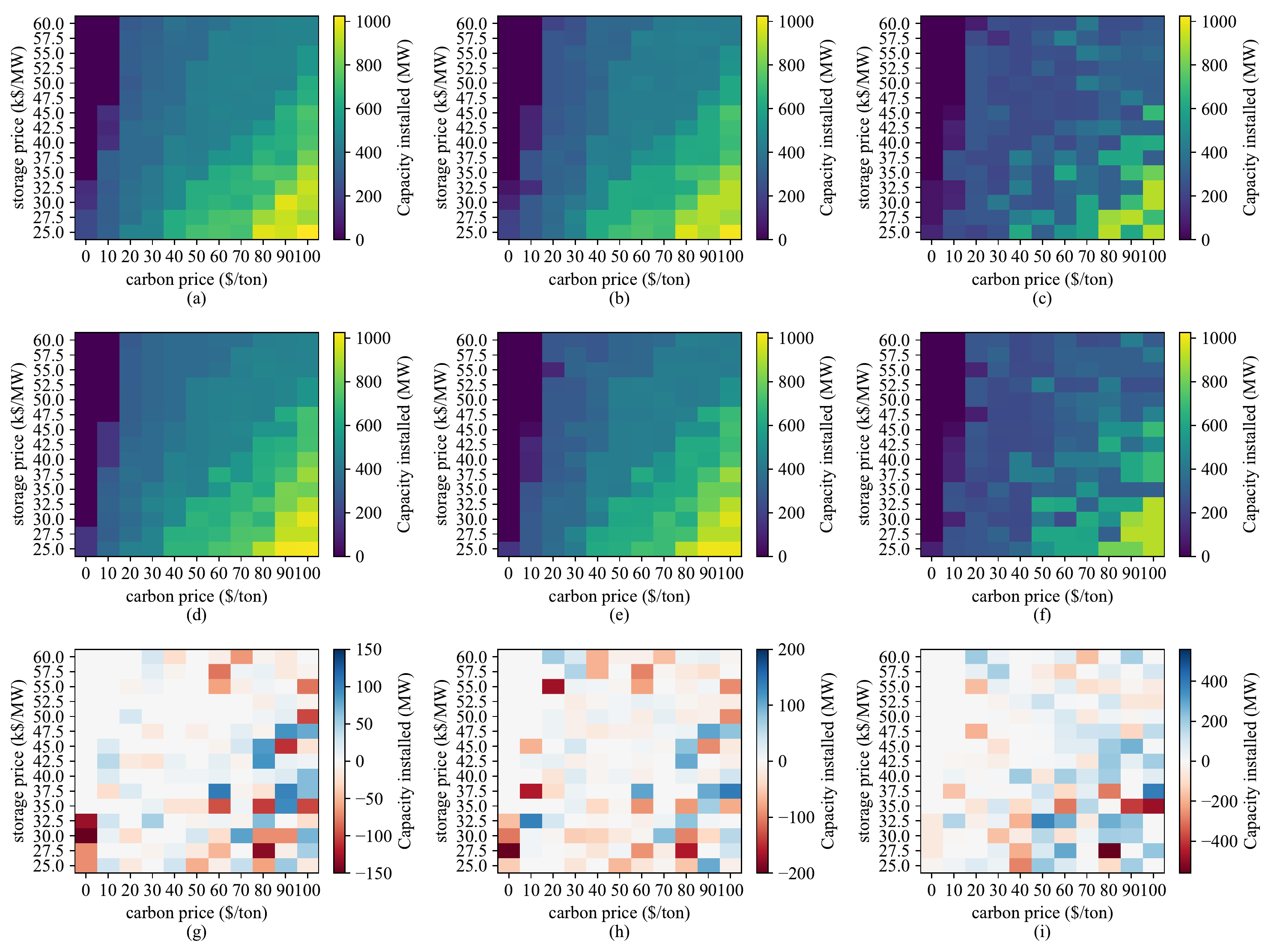}
    \caption{Optimal storage power as a function of investment perspective and emissions constraint. Rows: emissions unconstrained (a), (b), (c); emissions constrained (d), (e), (f); emissions constraint impact (g), (h), (i). Columns: vertically-integrated utility (a), (d), (g); philanthropic storage investor (b), (e), (h), profit-maximizing storage investor (c), (f), (i).}
    \label{3x3_storage}
\end{figure*}

\begin{figure*}[h]
    \centering
    \includegraphics[width=\linewidth]{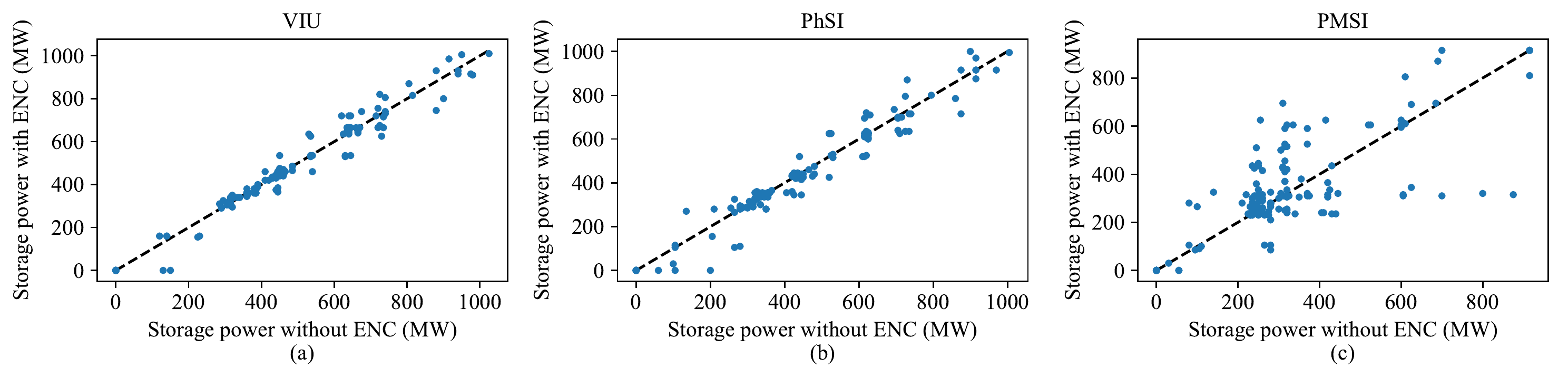}
    \caption{Optimal storage power with and without emissions-neutrality constraint, for a) vertically-integrated utility, b) philanthropic storage investor, and c) profit-maximizing storage investor.}
    \label{1x3_scatter}
\end{figure*}

\begin{figure*}[h]
    \centering
    \includegraphics[width=\linewidth]{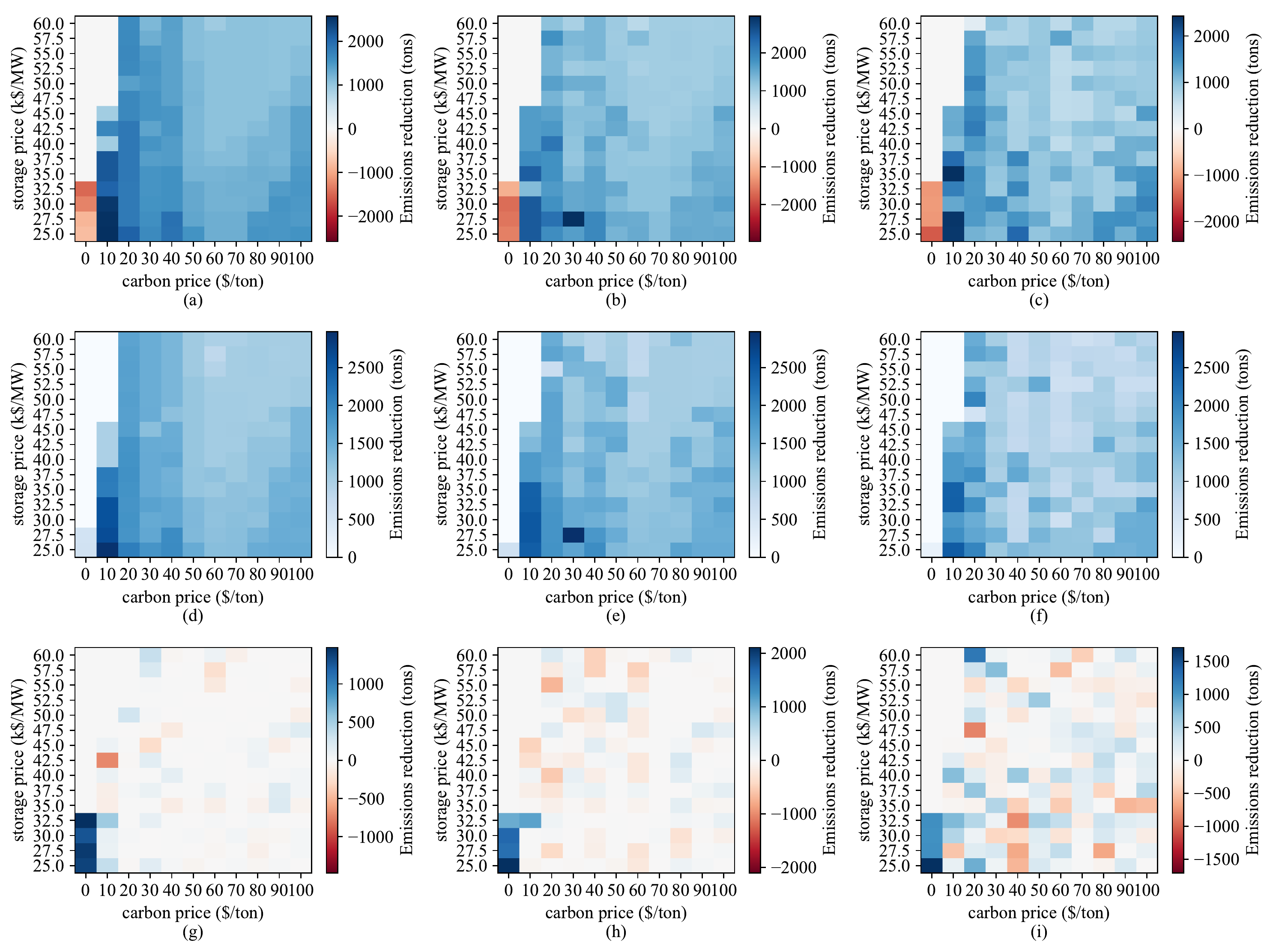}
    \caption{Emissions impact from storage as a function of investment perspective and emissions constraint. Rows: emissions unconstrained (a), (b), (c); emissions constrained (d), (e), (f); emissions constraint impact (g), (h), (i). Columns: vertically-integrated utility (a), (d), (g); philanthropic storage investor (b), (e), (h), profit-maximizing storage investor (c), (f), (i).}
    \label{3x3_emissions}
\end{figure*}

A case study is conducted using the Reliability Test System of the Grid Modernization Laboratory Consortium (RTS-GMLC) \cite{rts_gmlc}. The RTS-GMLC is based on the 1996 IEEE Reliability Test System, with an updated generation fleet, emissions curves for all fossil-fuel generators, the addition of renewable generation, and 365 days of hourly profiles for load and renewable generation. This system has a high penetration of renewables: over the course of the year, hydro represents 10.8\% of total energy demand, utility-scale solar PV 10.0\%, rooftop PV 5.7\%, and wind 19.0\%.

To obtain a `base case' system that is more representative of renewable penetrations in current high-renewables systems (or what will likely be encountered in the near future as penetration increases in lower-penetration systems), renewable sources were scaled down proportionally to a penetration of 30\%, measured on an annual energy basis. This renewable penetration is similar to that of the Irish grid at 33\% \cite{eirgrid2018}, notable for having the highest renewable penetration of large systems not within a larger synchronous AC network. The original RTS-GMLC renewable penetration (46\% by energy) is run as a `higher renewables' case.

For the planning problem, a set of five representative days is developed by using a k-means clustering algorithm on the 365 daily profiles, after their dimension was reduced by principal component analysis \cite{almaimouni2018}: 95\% of variance is captured via 13 principal components. When transforming the representative days back to full dimension, negative values of load or renewable generation are clipped to 0. Multi-period constraints in the system operator's problem (i.e. \eqref{binary_2}-\eqref{updown_2}, \eqref{soc_tracking}) are enforced cyclically; initial conditions are not assumed \textit{a priori}, but assumed to reflect the end of the current day \cite{olsen_2018tax}, and the net load in the first and last hours of the day are smoothed to avoid introducing unrealistic midnight ramps.

Dimensionality of the storage investment problem is reduced by fixing the storage power and energy ratios at a 4-hour duration, and by limiting storage investment to 10 candidate buses (out of the original 72); these candidate buses were selected by collecting a pool of candidate solutions to the VIU planning problem, and selecting all buses which had storage investment in any candidate solution.

To investigate the impact of price parameters on optimal storage quantities and the resulting emissions rates, a sensitivity analysis is conducted. Carbon prices of 0-100 \$/ton are tested in \$10/ton steps, and effective storage prices of 25,000-60,000 \$/MW-year are tested in \$2,500/MW-year steps. The effective storage price is an amortized value which depends on the capital cost of storage, equipment lifetime, financing rates, and presence of regulatory incentives.

Optimization problems were modeled using GAMS 25.0 and solved using CPLEX 12.8, on machines with at least 16 cores and 64 GB of RAM. Each storage siting and unit commitment problem was solved to an optimality gap of 0.1\% or better, and the heuristic algorithm optimizing both storage sizing and siting typically completed within 4 hours. As the heuristic features many independent storage siting optimizations, this time could further be reduced by running these optimizations in parallel.

\section{Results} \label{sec_results}

\subsection{Storage Quantity} \label{sec_res_stor}

Fig. \ref{3x3_storage} shows the optimal storage quantities found by a vertically integrated utility, a philanthropic storage investor, and a profit-maximizing storage investor, with sensitivity to both the effective price of energy storage (amoritized, per MW) and the price of carbon emissions. The first row shows the storage quantity without the emissions-neutrality constraint (ENC), the second row shows the storage quantity with the ENC, and the third row shows the difference from the ENC.

Decreasing storage prices and increasing carbon price quantities both tend to increase the optimal storage quantity, for all perspectives, both with and without an ENC. On average, the PhSI installs nearly as much storage as the VIU (93.3\% without the ENC, 91.8\% with the ENC), but the PMSI installs significantly less (67.5\% without the ENC, 70.7\% with the ENC). The effect of the ENC when there is no price on carbon is to reduce storage quantities in all three investment perspectives; on average the VIU selects 105 MW less storage, the PhSI selects 104 MW less, and the PMSI selects 45 MW less. However, the ENC also results in a significant decrease in emissions in these cases, as will be discussed in Section \ref{sec_res_emissions}.

The impact of the ENC on the optimal storage quantity selected by each type of investor is highly variable, especially in the PMSI perspective, as shown in Fig. \ref{1x3_scatter}. This variance is due to underlying nonconvexities of the optimization problem structures; all three investment perspectives feature discrete storage quantities and unit commitment variables, and solutions to these types of problems can be sensitive to changes in parameter values or new constraints. Additionally, since MIP models are typically solved to a specified optimality gap rather than global optimality, solvers may return solutions with similar objective function values but which differ sharply in their other decision variable values \cite{sioshansi2008}. Finally, the PMSI problem features arbitrage profit in its objective function, which is especially sensitive to changes in commitment solutions. The variance caused by these factors highlights the importance of conducting sensitivity analyses, to draw conclusions from repeated sampling rather than single point estimates.

Since the storage power solutions found in this case study are not normally distributed, non-parametric statistical methods must be used to draw conclusions in light of this high variance. The Wilcoxon signed-rank test is a non-parametric procedure which can be used to assess whether the means of a set of paired samples differ \cite{wilcoxon1945}; in this case, the installed storage quantities with and without the ENC are assessed to determine whether the null hypothesis that the means are the same can be rejected. Since there are several observations where the difference between the storage quantities is exactly zero, the procedure is adjusted as in \cite{pratt1959}.

When there is a nonzero price on carbon, there is no statistically significant difference in storage quantity for either the VIU or PhSI perspective ($p=0.90$ and $p=0.10$, respectively), while for the PMSI perspective there is a statistically significant but modest increase in the optimal storage quantity ($p=0.0095$, average increase of 17.5 MW). This seemingly counter-intuitive result is due to the fact that storage dispatch and generator commitment/dispatch are considered simultaneously, so the emissions constraint applies to the unit commitment variables which, once fixed, define the time- and location-varying generator supply curves. The optimal storage quantity is ultimately a function of these supply curves, and in these case studies the ENC seems to result in supply curves which offer maximum PMSI profitability at greater storage quantities.

Overall, the effect of the ENC on storage quantity appears to be much weaker than the effects of the storage price, the carbon price, and the variance associated with re-solving mixed-integer linear programs with perturbed parameters.

\subsection{Emissions Quantity} \label{sec_res_emissions}

Fig. \ref{3x3_emissions} shows the emissions impact of the storage quantities shown in Fig. \ref{3x3_storage}. Although the emissions impact of storage is seen to be beneficial (reducing annual emissions) when there is a carbon price, in the absence of a carbon price the impact of adding storage is to increase emissions, as shown in Figs. \ref{3x3_emissions}(a), \ref{3x3_emissions}(b), and \ref{3x3_emissions}(c). Emissions are increased by an average of 2.3\% for the VIU perspective, 3.1\% for the PhSI perspective, and 2.5\% for the PMSI perspective. On average, renewable spillage is reduced by 81\%, but this beneficial effect is outweighed by a 7\% increase in energy generated from coal, coupled with a 10\% decrease in energy generated from natural gas.

The impact of the ENC is to significantly reduce the system emissions in cases without a carbon price (on average 2.9\% for the VIU, 3.4\% for the PhSI, 2.6\% for the PMSI), to below the emissions baseline (no-storage) while having little effect on the emissions in cases with a carbon price, as shown in Figs. \ref{3x3_emissions}(g), \ref{3x3_emissions}(h), and \ref{3x3_emissions}(i) (The only perspective which saw a statistically significant result in the presence of a carbon price was the PMSI, with an average emissions reduction of 0.2\%). Emissions reductions in the cases without a carbon price were achieved with only minor cost increases: on average 0.11\% for the VIU, 0.09\% for the PhSI, 0.05\% for PMSI, and no statistically significant change in the priced GHG cases ($p>0.9$ in the Wilcoxon test).

Although seemingly counter-intuitive, enforcement of the ENC on a daily basis can sometimes tend to increase overall emissions, as seen in Figs. \ref{3x3_emissions}(g), \ref{3x3_emissions}(h), and \ref{3x3_emissions}(i). In the absence of an ENC, there may be certain days of the year where storage tends to increase emissions and certain days where it tends to decrease emissions, based on the particular profiles of electricity demand and renewable generation. When the ENC binds on the emissions-increasing days, and lowers the economic value of energy storage, the optimal quantity of storage to invest in may change. The net impact of this change in storage quantity may be negative if the beneficial effect on the ENC in emissions-increasing days is outweighed by the detrimental effect on emissions-decreasing days (i.e., the new storage quantity is not as effective at decreasing emissions on those days).

The quality of the PhSI solutions can be evaluated by investigating the lower bound of the two relaxations described in Section \ref{sec_solution_quality}. The solutions found by the heuristic in Algorithm \ref{alg_1} are typically within 0.1\% of the best lower bound found by the VIU relaxation. Although the PCSLE formulation is a tighter relaxation of the PhSI problem, in practice the progress of the best lower-bound is relatively slow; the best lower bound found by the PCSLE formulation after 24 hours is still below the best lower bound of the VIU, which typically solves to 0 optimality gap within one hour.

\subsection{Higher Renewables Case}

When the methods are applied to the original RTS-GMLC test system with 46\% penetration of renewables by energy, the results are broadly similar, but with a few differences. In the higher renewables case the addition of energy storage reduces the average system emissions, even without a price on carbon: by 1.0\% in the VIU perspective, 0.5\% in the PhSI perspective, and by 0.4\% in the PMSI perspective. This is consistent with the results from \cite{craig2018}, which shows that the negative emissions effects from storage become positive as power systems increase their penetrations of low- and no-carbon generation sources.

However, there are still some days for which storage increases emissions, so the addition of the ENC still has an effect on the operation of these days and the optimal storage quantity: the ENC reduces emissions by 2.6\% for the VIU, 3.1\% for the PhSI, and 3.0\% for the PMSI. Notably, these reductions are similar in magnitude to the emissions reductions observed in the case with 30\% renewables, implying that although the average effect of storage on emissions may be better in higher-renewables cases, the beneficial impact of the ENC may remain relatively constant. The increase in total cost is modest: 0.12\% in the VIU case, 0.14\% in the PhSI case, and 0.19\% in the PMSI case. When there is a price on carbon, there is no statistically significant difference in cost or emissions in any of the investment perspectives ($p>0.4$).

In all investment perspectives, the optimization prioritizes installing storage first at bus 309 (downstream of congested line C6), due to the economic opportunities that this congestion presents. This is contrasted with the 30\% penetration `base case', which does not exhibit congestion in any of the five representative days (although it is present at this line on some days when the number of representative days is expanded from five to ten).

\section{Conclusion} \label{sec_conclusion}

The addition of energy storage to power systems operating without a price on carbon can increase system emissions, even in systems with high penetrations of renewables (30\% in this case study). Even in a system with much higher penetration (46\% in the `higher renewables' case), the impact may be on average only minorly beneficial, as storage tends to increase emissions in some days and decrease it in others. This effect persists over a wide range of energy storage prices, suggesting that policies which subsidize energy storage installation may not result in lowered emissions, absent complementary policies governing storage operation.

If pricing GHG emissions is not feasible (\textit{e.g.} due to political obstacles), then adding an emissions-neutrality constraint to the operation of energy storage is shown to have a significant beneficial effect on system emissions (approximately 3\%, depending on the storage investment perspective), with only a minor impact on overall cost (approximately 0.1-0.2\%, depending on the storage investment perspective and the penetration of renewables). A `philanthropic' (profit-constrained) storage investor tends to invest in significantly more energy storage than a profit-maximizing storage investor, illustrating the benefits that a socially-minded storage investment entity can provide, even in the presence of profitability and emissions-neutrality constraints.

While these effects have been shown in this particular case study, future work can examine the sensitivity of these results to more general cases. Factors which may contribute to different outcomes in terms of storage investment and emissions may include:

\begin{itemize}
    \item differences in market design, such as participation of energy storage in the reserves market(s),
    \item varying total penetrations of non-emitting resources, or varying ratios of solar, wind, and nuclear generation, and
    \item consideration of different storage durations, or of investment in a mixture of storage with varying durations.
\end{itemize}

\bibliographystyle{IEEEtran}
\bibliography{bibliography}

\section*{Appendix A: Dual problem of Transmission Constrained Economic Dispatch}

\subsection{Creating a Profit-Constrained Single-Level Equivalent} \label{sol_quality_sle}

Though the PhSI objective \eqref{strat_storage_phil_obj} and operational constraints \eqref{binary_2}-\eqref{emissions_constraint} can be formulated as a single-level MILP problem, the profit definition \eqref{profit_def} contains both primal and dual decision variables; therefore, these variables must be optimized simultaneously.
The system operator's unit commitment problem is necessarily non-convex due to binary commitment variables; strong duality does not in general hold for non-convex problems. However, if the binary variables and constraints are moved to the upper-level of the bi-level formulations described in Section \ref{sec_solution_method}, the lower level becomes a transmission-constrained economic dispatch (TCED) problem, which can be represented linearly using the DC power flow approximation.
A linear (and therefore convex) lower-level problem can be replaced by a series of constraints in an upper-level problem, creating a single-level equivalent of the profit-constrained storage investment problem \cite{dvorkin2017}.

Moving the binary variables and constraints from the lower-level to the upper-level represents a relaxation of the original bi-level problem, since the values of the binary variables are no longer constrained by the lower-level unit commitment problem, and no additional constraints are introduced. Therefore, the best lower bound on the relaxed problem provides a lower bound for the best value of the original PhSI problem.

Since the original TCED problem is linear, strong duality is guaranteed, so the value of the primal and dual objective functions are equal for a set of primal and dual variables ($\Omega^\lambda$) that are primal and dual optimal, respectively. The dual problem of the TCED is given in \eqref{dual_objective}-\eqref{dual_sren_cons}, with the primal variables to which each dual constraint corresponds listed in parentheses. The strong duality constraint for the TCED problem is given in \eqref{strong_duality}. Although this single-level equivalent contains bi-linear terms, Section A.2 demonstrates how a MILP approximation can be obtained, which can provide a lower-bound on the optimum.

\begin{gather}
    \max_{\Omega^{\text{D}}} \hat{C}^{\text{dual}} := \alpha \left [ \sum_{t \in T} \sum_{i \in I} E_{i}^{\text{min}} u_{i,t} + E_{i}^{\text{su}} v_{i,t} - \chi E^{\text{baseline}} \right ]
    \nonumber \\
    - \sum_{l \in L} \sum_{t \in T} \left [
    f_{l}^{\text{max}} \left ( \underline{\gamma}_{l,t} + \overline{\gamma}_{l,t} \right )
    \right ]
    - \sum_{i \in I} \sum_{s \in S} \sum_{t \in T}
    g_{i,s}^{\text{max}} u_{i,t}
    \overline{\delta}_{i,s,t}
    \nonumber \\
    + \sum_{b \in B} \sum_{t \in T} \Bigg [
    \nonumber
    \lambda_{b,t}  \left ( d_{b,t} - w_{b,t}
    - \sum_{i \in I} m_{i,b}^{\text{unit}} g_{i}^{\text{min}} u_{i,t} \right )
    \nonumber \\
    - 
    Q_{b}^{\text{max}}
    \overline{\xi}_{b,t}
    -
    J_{b}^{\text{max}}
    \left ( \overline{\rho}_{b,t}^{\text{chg}} + \overline{\rho}_{b,t}^{\text{dis}} \right )
    -
    w_{b,t}
    \overline{\phi}_{b,t}
    \Bigg ]
    \label{dual_objective} \displaybreak[0] \\
    \text{subject to:} \nonumber \\
    \alpha h_{i,s} + b_{i,s} - \sum_{b \in B}
    \left [ m_{i,b}^{\text{unit}} \lambda_{b,t} \right ] - \underline{\delta}_{i,s,t} + \overline{\delta}_{i,s,t}
    = 0
    \hspace{8pt} \left ( g_{i,s,t} \right )
    \label{dual_gist_cons} \\
    \frac{1}{\eta^{\text{dis}}} \kappa_{b,t} + \overline{\rho}_{b,t}^{\text{dis}} - \underline{\rho}_{b,t}^{\text{dis}} - \lambda_{b,t} = 0
    \hspace{8pt} \left ( J_{b,t}^{\text{dis}} \right )
    \label{dual_jdis_cons} \\
    \lambda_{b,t} + \overline{\rho}_{b,t}^{\text{chg}} - \underline{\rho}_{b,t}^{\text{chg}} - \eta^{\text{chg}} \kappa_{b,t} = 0
    \hspace{8pt} \left ( J_{b,t}^{\text{chg}} \right )
    \label{dual_jchg_cons} \\
    \sum_{l \in L} y_{l} m_{l,b}^{\text{line}} \beta_{l,t} = 0
    \hspace{8pt} \left ( \theta_{b,t} \right )
    \label{dual_theta_cons2} \\
    \sum_{b \in B} \left ( m_{l,b}^{\text{line}} \lambda_{b,t} \right ) + \beta_{l,t} - \underline{\gamma}_{l,t} + \overline{\gamma}_{l,t} = 0
    \hspace{8pt}
    \left ( f_{l,t} \right )
    \label{dual_f_cons} \\
    \kappa_{b,t} - \kappa_{b,t+1} + \overline{\xi}_{b,t} - \underline{\xi}_{b,t}
    = 0
    \hspace{8pt} \left ( Q_{b,t} \right )
    \label{dual_q_cons} \\
    \lambda_{b,t} + \overline{\phi}_{b,t} -  \underline{\phi}_{b,t}
    = 0
    \hspace{8pt} \left ( s_{b,t}^{\text{ren}} \right )
    \label{dual_sren_cons} \\
    \alpha, \gamma, \delta, \zeta, \xi, \rho, \phi \ge 0 \nonumber
\end{gather}

\begin{equation}
    \hat{C}^{\text{dual}} = \hat{C}^{\text{primal}} := \sum_{i \in I} \sum_{s \in S} \sum_{t \in T} b_{i,s} g_{i,s,t}
    \label{strong_duality}
\end{equation}

Therefore, the profit-constrained single-level equivalent of the relaxed PhSI problem is given by the original objective function \eqref{strat_storage_phil_obj_again}, the upper-level profit constraint \eqref{profit_constraint_again2}, unit commitment constraints \eqref{commitment_constraints_again2}, and lower-level constraints \eqref{primal_feasibility_constraints_again2}-\eqref{dual_feasibility_constraints_again2}.

\begin{gather}
    \min_{\Omega^{\text{C}}, \Omega^{\text{D}}, \Omega^{\text{S}}, \Omega^{\lambda}} C^{\text{batt}} + \sum_{a \in A} C_{a}^{\text{gen}}
    \label{strat_storage_phil_obj_again} \\
    \text{subject to:} \nonumber \\
    \text{Equation \eqref{strat_storage_phil_profit_cons}}
    \label{profit_constraint_again2} \\
    \text{Equations \eqref{binary_2}-\eqref{updown_2}}
    \hspace{20pt} \forall a \in A
    \label{commitment_constraints_again2} \\
    \text{Equations \eqref{seg_limit}-\eqref{emissions_constraint}}
    \hspace{20pt} \forall a \in A
    \label{primal_feasibility_constraints_again2} \\
    \text{Equation \eqref{strong_duality}}
    \hspace{20pt} \forall a \in A
    \label{strong_duality_again2} \\
    \text{Equations \eqref{dual_gist_cons}-\eqref{dual_sren_cons}}
    \hspace{20pt} \forall a \in A
    \label{dual_feasibility_constraints_again2}
\end{gather}

\subsection{Creating a MILP Approximation} \label{sol_method_milp}

The profit-constrained single-level equivalent presented in \eqref{strat_storage_phil_obj_again}-\eqref{dual_feasibility_constraints_again2} contains several non-linear terms. The profit definition \eqref{profit_def} contains the product of continuous lower-level dual variables ($\boldsymbol{\lambda}$) and primal variables ($\boldsymbol{J}^{\text{chg}}, \boldsymbol{J}^{\text{dis}}$), while the strong duality constraint contains products of continuous lower-level dual variables and upper-level binary and continuous variables. First, the profit constraint can be converted from a product of lower-level primal and dual variables to a product of lower-level dual and upper-level variables using complementary slackness conditions, as shown in \eqref{profit_de_lambda} \cite{dvorkin2017}.

\begin{gather}
    \sum_{b \in B} \sum_{t \in T} \lambda_{b,t,a} \left ( J_{b,t,a}^{\text{dis}} - J_{b,t,a}^{\text{chg}} \right )
    =
    \hspace{100pt} \nonumber \\ \hspace{50pt}
    \sum_{b \in B} \sum_{t \in T} \Big [
    Q_{b}^{\text{max}}
    \overline{\xi}_{b,t}
    +
    J_{b}^{\text{max}}
    \left ( \overline{\rho}_{b,t}^{\text{dis}} + \overline{\rho}_{b,t}^{\text{chg}} \right )
    \Big ]
    \label{profit_de_lambda}
\end{gather}

Next, the continuous upper-level variables $Q_{b}^{\text{max}}$ and $J_{b}^{\text{max}}$ can be approximated by integer variables representing discrete storage quantities. These integer variables can be equivalently represented by a summation of binary variables in order to convert the products of integer and continuous variables to the products of binary and continuous variables. A binary expansion is used in \eqref{discretize_q}-\eqref{discretize_j} to reduce dimensionality as compared to a unary expansion. Although the representation of integer variables by a binary expansion is not in general more efficient than an integer representation \cite{owen2002}, this allows the use of the big-M method, and has been shown to be more effective than non-linear solvers or the use of McCormick envelopes in solving bi-linear problems containing the product of continuous and integer variables in constraints \cite{gupte2013}.

\begin{gather}
    Q_{b}^{\text{max}} \approx \sum_{n} 2^{n} x_{b,n}^{a} \Delta Q
    \label{discretize_q} \\
    J_{b}^{\text{max}} \approx \sum_{n} 2^{n} x_{b,n}^{b} \Delta J
    \label{discretize_j}
\end{gather}

After discretization and binary expansions of the storage variables, the only non-linear terms remaining are products of continuous lower-level dual variables and upper-level binary variables. These products are linearized using the big-M approximation method, and the problem reduces to a MILP approximation of the original non-linear problem.

\begin{IEEEbiography}[{\includegraphics[width=1in,height=1.25in,clip,keepaspectratio]{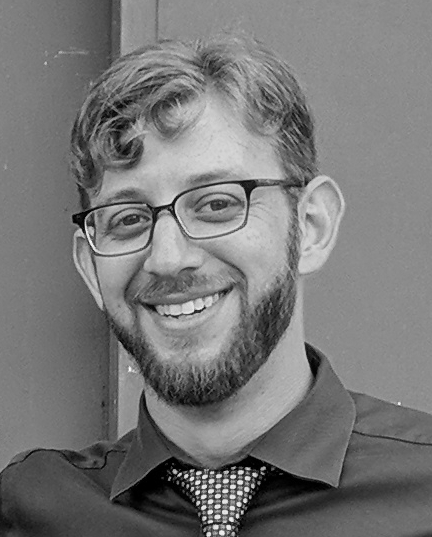}}]{Daniel Olsen} (S'14-M'19) is a Research Scientist at Intellectual Ventures. His research interests include policies and planning for low-carbon power systems, integration strategies for high penetrations of intermittent renewable resources, and the valuation of distributed flexibility resources.

He holds a Ph.D. from the University of Washington and a B.Sc. in Mechanical Engineering and Electric Power Engineering from Rensselaer Polytechnic Institute.
\end{IEEEbiography}

\begin{IEEEbiography}[{\includegraphics[width=1in,height=1.25in,clip,keepaspectratio]{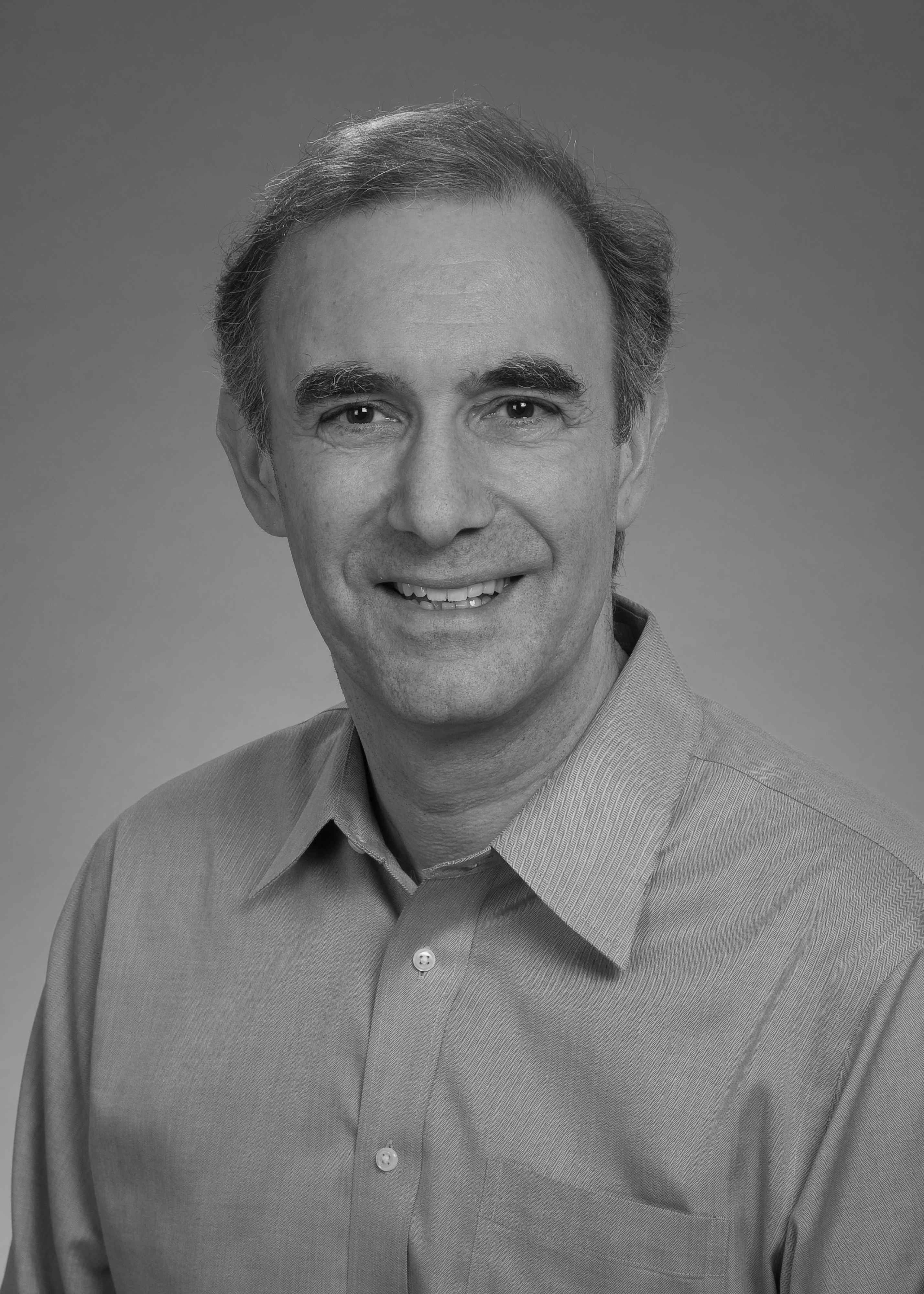}}]{Daniel Kirschen} (M’86–SM’91–F’07) 
is the Donald W. and Ruth Mary Close Professor of Electrical Engineering at the University of Washington. His research focuses on the integration of renewable energy sources in the grid, power system economics and power system resilience. Prior to joining the University of Washington, he taught for 16 years at The University of Manchester (UK). Before becoming an academic, he worked for Control Data and Siemens on the development of application software for utility control centers. He holds a PhD from the University of Wisconsin-Madison and an Electro-Mechanical Engineering degree from the Free University of Brussels (Belgium). He is the author of two books and a Fellow of the IEEE.
\end{IEEEbiography}

\end{document}